\documentclass[a4ja,10pt]{article}
\newtheorem{thm}{Theorem}[section]
\newtheorem{prop}[thm]{Proposition}

\newtheorem{pf}{Proof}

\newtheorem{defn}[thm]{Definition}

\usepackage{amssymb}
\usepackage{latexsym}

\newcommand{\cpn}{{{\mathbb C}{\mathbb P}^{n}}}
\newcommand{\sts}{{\mathbb C}{\mathbb P}^{3}}

\newcommand{\proj}{{\rm Proj}(S)} 

\newcommand{\oproj}{{\mathcal O}_{{\rm Proj}(S)}} 
 
\begin{document}


\begin{center}
{\Large\bf Symbol calculus on a projective space}
\end{center}

\par\bigskip
\begin{center}
Naoya MIYAZAKI\\
Department of Mathematics, 
Keio University, \\
Yokohama, 223-8521, JAPAN
\end{center}



\par\medskip\noindent
{\small {\bf Abstract}:
In this article, we introduce symbol calculus on a projective scheme. 
Using holomorphic Poisson structures, 
we construct deformations of ring structures  
for structure sheaves on projective spaces $\cpn$'s. 

\par\medskip\noindent
\par\medskip
\noindent
{\bf Acknowledgements:} 
The author would like to dedicate the present article 
to Professor Akira Yoshioka on his 60th birthday.
This research is partially supported by 
JSPS Grant-in-Aid for Scientific Research, 
and the academic fund of Keio University. 

\section{Introduction}
The terminology ``symbol calculus" is used in Fourier analysis 
and pseudo-differential operator 
in order to study partial differential equations. 
Especially, symbol calculus of pseudo-differential operator 
with respect to elliptic operator gives fruitful contribution 
to the index theorem, i.e. 
an important role is played 
by pseudo-differential operator 
in the index theory for elliptic operators appearing geometry, 
where pseudo-differential 
operators are used to extend the class of possible deformation 
of an elliptic operator which has essential topological datas 
of the base manifold. 

In this article, 
using holomorphic Poisson structures, 
we consider symbol calculus on a projective scheme. 
Especially, we give concrete examples using $\cpn$. 
As mentioned above, structure sheaves on algebraic varieties have 
important and essential feature which plays 
crucial role to analyze their fundamental properties with respect to 
the base varieties. 
Using a holomorphic Poisson structure on the base variety, 
we construct symbol calculus on the structure sheaf. 
We here state the main theorems of the present article. 

\begin{thm}\label{thm1}
Assume that $Z=[z_0:z_1:\ldots:z_n]$ is 
the homogeneous coordinate system of $\cpn$, 
and $\Lambda=\sum_{\alpha,\beta}\stackrel{\leftarrow}{\partial_{Z_{{\alpha}}}}
\Lambda^{{\alpha},{\beta}} 
        \stackrel{\rightarrow}{\partial_{Z_{{\beta}}}}$
defines a holomorphic skew-symmetric biderivation
\footnote{We use Einstein's convention unless confusing.}  
of order zero acting on the structure sheaf ${\mathcal O}_\cpn$ 
satisfying the Jacobi rule
\footnote{The biderivation used here might be called 
an ``algebraic Poisson structure." I am not sure. } 
and an assumption below: 
{ 
\begin{eqnarray}\label{lambda-relation}
\nonumber
&&
\bigr(\stackrel{\leftarrow}{\partial_{Z_{{\alpha}_1}}}
\Lambda^{{\alpha}_1,{\beta}_1} 
        \stackrel{\rightarrow}{\partial_{Z_{{\beta}_1}}}\bigr)
\cdots
\bigr(\stackrel{\leftarrow}{\partial_{Z_{{\alpha}_k}}}
\Lambda^{{\alpha}_k,{\beta}_k} 
        \stackrel{\rightarrow}{\partial_{Z_{{\beta}_k}}}\bigr)\\
&=&\stackrel{\longleftarrow}{\partial_{Z_{{\alpha}_1\ldots{\alpha}_k}}}
\Lambda^{{\alpha}_1, {\beta}_1}\cdots
\Lambda^{{\alpha}_k,{\beta}_k} 
        \stackrel{\longrightarrow}{\partial_{Z_{{\beta}_1\ldots{\beta}_k}}}, 
\end{eqnarray}
}
where 
$$
f(Z)\stackrel{\longleftarrow}{\partial_{Z_{{\alpha}_1\ldots{\alpha}_k}}}
\quad
(\mbox{resp}. ~ 
\stackrel{\longrightarrow}{\partial_{Z_{{\beta}_1\ldots{\beta}_k}}} 
g(Z))
$$  
means 
$$\partial_{Z_{\alpha_1}}\partial_{Z_{\alpha_2}}
 \cdots \partial_{Z_{\alpha_k}}f(Z)
\quad(\mbox{resp.} ~ 
\partial_{Z_{\beta_1}}\partial_{Z_{\beta_2}}
\cdots \partial_{Z_{\beta_k}}g(Z)).$$ 
Then, for any point ``$p$" and germs $[f(Z)]_{~{}_{\tilde{}}}$, 
$~[g(Z)]_{~{}_{\tilde{}}}$ of 
the stalk ${\mathcal O}_{\cpn,p}[[\natural]]$
\footnote{Here $[[\natural]]$ denotes  
either $[\mu,\mu^{-1}]]$ or $[\mu,\mu^{-1}]$, 
and we have to choose carefully in context. }, 
\begin{equation}\label{sharp-product}
\begin{array}{ll}
& f(Z)\# g(Z) \\
:=&\sum_{k=0}^\infty \frac{1}{k!}\left(\frac{\mu}{2}\right)^k
\Lambda^{\alpha_1\beta_1}
\Lambda^{\alpha_2\beta_2}
\cdots \Lambda^{\alpha_k\beta_k}
\partial_{Z_{\alpha_1}}\partial_{Z_{\alpha_2}}\cdots \partial_{Z_{\alpha_k}}f(Z)\partial_{Z_{\beta_1}}\partial_{Z_{\beta_2}}\cdots \partial_{Z_{\beta_k}}g(Z)
\end{array}
\end{equation}
defines a non-commutative and associative ring structure, 
where $\mu$ is a formal parameter. 
\end{thm}
We also have 
\begin{thm}\label{thm1-5}
Under the same assumptions and notations of Theorem~\ref{thm1}, 
the product $\#$ induces globally defined non-commutative, associative 
product on the sheaf-cohomology space  
$\sum_{k=0}^\infty H^0(\cpn,{\mathcal O}_\cpn(k))[[\natural]]$ 
where $\mu$ can be specialized a scalar (for example $\mu=1$). 
\end{thm}
Using the product $\#$, we also have 
\begin{thm}\label{thm2} 
Suppose that the same assumptions of Theorem~\ref{thm1}. 
Let $A[Z]:=ZA{}^tZ$ be a quadratic form with homogeneous degree 2. 
Then we have 
\begin{equation}\label{star-exponential}
e_{\#}^{\frac{1}{\mu}A[Z]}=\det{}^{-1/2}
\left(\frac{e^{\Lambda A}+e^{-\Lambda A}}{2}\right)\cdot 
e^{1/\mu\left(\frac{\Lambda^{-1}}{\sqrt{-1}}\tan (\sqrt{-1}\Lambda A)\right)[Z]}\in \sum_{k=0}^\infty H^0(\cpn,{\mathcal O}_\cpn(k)[\mu,\mu^{-1}]]).
\end{equation}
\end{thm}

\section{Outlines of Proofs of 
Theorems \ref{thm1}, \ref{thm1-5} and \ref{thm2}}
In this section we give outline of proofs 
of our main theorems (cf. \cite{miyazaki-preparation}). 
Before giving proofs, 
we recall fundamentals of projective scheme.  
In fact, we need a slight modification of the standard theory of scheme. 

Let $S=\oplus_{n=0}^\infty S_n$ be a graded commutative ring. 
Then, $S_0$ is obviously commutative and $S$ is an $S_0$-algebra. 
It is well-known that the homogeneous ideal $S_+:=\oplus_{n=1}^\infty S_n$  
is called the {\it irrelevant ideal}. And the following is well-known: 
\begin{prop}\label{noether}
A graded commutative ring $S$ is noetherian if and only if 
$S_0$ is noetherian and $S$ is finitely generated by $S_1$ as an $S_0$-algebra. \end{prop}
It is also well-known that a projective shceme 
$$\proj :=\{\frak{p}:\mbox{ a homogeneous prime ideal}~|~ 
\neg(S_+\subset \frak{p} ) \} $$
admits the canonical scheme structure in the following way: 
Set 
$$D_+(f):=\bigr\{\frak{p} \in \proj ~|~ \neg(f\in \frak{p} ) \bigr\}, $$ 
for any homogeneous element $f\in S_d$ with degree $d$, 
then the family $\{D_+(f)\}_{f\in S_d,~d\in {\mathbb Z}_{\geq 0 }}$ forms 
a basis of open sets. Hence it gives the canonical topology $\frak{O}_{\proj}$ 
(that is, the Zariski topology) for $\proj$. 
Note that 
$\neg(f\in \frak{p} )$ means $f(\mbox{point}_{\frak{p}})\not= 0$, intuitively. 
We also set 
\begin{equation}\label{5.7}
\begin{array}{ll}
\Gamma(D_+(f),{\oproj}):=\bigr\{g/f^m~|~g\in S_m,~m\geq 0\bigr\},\\
\oproj:\frak{O}_{\proj} \ni D_+(f)  
\mapsto \Gamma(D_+(f),\oproj)\in {\rm\mathbf Mod}. 
\end{array}
\end{equation}
The functor above is well-known as the {\it structure sheaf}. 
We remark that when $g\in S_0$, we easily see that $g/1=fg/f~(fg\in S_d)$. 
Hence we may consider $g/f^m~(m\geq 1 )$ instead of $g/f^m~(m\geq 0 )$. 
We obtain that R.H.S. of (\ref{5.7}) is a part of degree $0$ of localization 
$S_f$ of $S$ by a product closed set $\{f^{\ell}\}_{\ell=0,1,2,\ldots } $. 
We denote it by $(S_f)_0$ or $S_{(f)}$. Strictly speaking,  
for any homogeneous element $f$ with ${\rm deg}(f)=d$, 
$$(S_f)_0:=S_{(f)}:=\bigr\{g/f^m~|~g\in S_{md},~m\geq 0 \bigr\}.$$
Hence we proved that 
\begin{prop}\label{scheme}
As for $D_+(f)$, 
\begin{equation}\label{affine-scheme}
(D_+(f),\oproj|_{D_+(f)})\cong {\rm Spec}(S_{(f)}).
\end{equation} 
Thus, $\proj$ is obtained by glueing of affine schemes. 
It indicates that $(\proj,\oproj)$ is a scheme in the genuin sense. 
\end{prop}
Next we consider cohomology of quasi-coherent sheaf over $\proj$. 
Assume that a graded ring $S$ is generated by $S_1$ as an $S_0$-algebra. 
For instance 
$$S=R[z_0,z_1,\ldots,z_n],~S_0=R,~S_1=\rm\{a\in S~|~{\rm deg}(a)=1\rm\}.$$
As for a quasi-coherent sheaf 
\footnote{A sheaf ${\mathcal F}$ is quasi-coherent if and only if there is a
pre-sheaf exact sequence ${\mathcal O}_U^{\oplus I}\to {\mathcal O}_U^{\oplus J}\to {\mathcal F}\to 0 $.
A sheaf ${\mathcal F}$ is coherent if and only if there is a
pre-sheaf exact sequence ${\mathcal O}_U^{\oplus n}\to {\mathcal F}\to 0~~
(n\in {\mathbb N}) $.}
${\mathcal F}$, we set 
$${\mathcal F}(m)[[\natural]]:={\mathcal F}\otimes_{\oproj} \oproj(m)[[\natural]],$$ 
and define 
\begin{equation}\label{5.21}
\Gamma_*({\mathcal F}):=
\oplus_{m\in {\mathbb Z}} \Gamma(X,{\mathcal F}(m))[[\natural]],~~
{\rm deg}(a):=m,~~(\forall a \in \Gamma(X,{\mathcal F}(m)). 
\end{equation}
Then we see that $\Gamma_*({\mathcal F}[[\natural]])$ 
is a graded $\Gamma(\oproj [[\natural]])$-module. 
For any element $f\in S_d$, we set $\alpha_d(f):=a/1$. 
Then it is well-known that 
\begin{prop}\label{thm-alpha}
The map $\alpha_d$ obtained above defines a homomorphism 
\begin{equation}\label{alpha}
\alpha_d(f):
S_d[[\natural]]\ni a \mapsto a/1 \in S(d)_{(f)}=\Gamma(D_+(f),\oproj (d) [[\natural]]).
\end{equation}
A family $\{\alpha_d(f)\}_{f\mbox{:homogeneous}}$ induces 
a module homomorphism 
\begin{equation}\label{5.22-1}
\alpha_d~:S_d[[\natural]]\to\Gamma(\proj,\oproj(d) [[\natural]]). 
\end{equation}
Hence, using the module homomorphisms $\{\alpha_d\}$, 
a graded ring homomorphism 
\begin{equation}\label{5.22}
\alpha:=\oplus_{n=0}^\infty 
\alpha_d~:S=\oplus_{n=0}^\infty S_d[[\natural]]\to\Gamma(\proj,\oproj [[\natural]])  
\end{equation}
can be defined for any quasi-coherent sheaf ${\mathcal F}$. 
Thus, $\Gamma(\proj,\oproj [[\natural]])$ admits 
a graded $S[[\natural]]$-module structure. 
\end{prop}
\begin{defn}
We denote the pair $({\rm Proj}(S),{\mathcal O}_{{\rm Proj}(S)}[[\natural]])$ 
by ${\rm Proj}(S)[[\natural]]$. 
\end{defn}
We are also interested in $\Gamma_*({\mathcal F}[[\natural]])_{(f)}$. As for any element $f\in S_d$, and $x\in\Gamma(\proj,{\mathcal F}[[\natural]])$, we see 
$x/f^n\in \Gamma_*(\proj,{\mathcal F}[[\natural]])_{(f)}.$  
We denote the restriction of $x$ to $D_+(f)$ by $x|_{D_+(f)}$. 
Then matching the degrees of $x|_{D_+(f)}$ and $(a_d(f)|_{D_+(f)})^n$ 
we get $x|_{D_+(f)}/\bigr(a_d(f)|_{D_+(f)})\bigr)^n.$ 
\begin{equation}\label{5.23-1}
\beta_{(f)}:\Gamma_*(\proj,{\mathcal F}[[\natural]])_{(f)}\ni \frac{x}{f^m}
\mapsto 
\frac{x|_{D_+(f)}}{(\alpha(f)|_{D_+(f)})^m}
\in \Gamma(D_+(f),{\mathcal F}[[\natural]]). 
\end{equation}
As similarly for $\alpha$, for any homogeneous element $g\in S_e$, 
we obtain a diagram:
\par\bigskip
$\qquad\qquad$
\unitlength 0.1in
\begin{picture}( 29.9500, 19.8500)(  4.0000,-21.9500)
\put(4.0000,-4.0000){\makebox(0,0)[lb]{$\!\!\!\!\!\!\!\!\!\!\!\!\!\!\!
\Gamma_*({\mathcal F}[[\natural]])_{(f)}$}}%
%
\special{pn 8}%
\special{pa 1000 350}%
\special{pa 3000 350}%
\special{fp}%
\special{sh 1}%
\special{pa 3000 350}%
\special{pa 2934 330}%
\special{pa 2948 350}%
\special{pa 2934 370}%
\special{pa 3000 350}%
\special{fp}%
\put(32.0000,-3.8000){\makebox(0,0)[lb]{$\Gamma(D_+(f), {\mathcal F}[[\natural]])$}}%
%
\special{pn 8}%
\special{pa 600 590}%
\special{pa 600 1990}%
\special{fp}%
\special{sh 1}%
\special{pa 600 1990}%
\special{pa 620 1924}%
\special{pa 600 1938}%
\special{pa 580 1924}%
\special{pa 600 1990}%
\special{fp}%
\put(4.0000,-22.4000){\makebox(0,0)[lb]{$\!\!\!\!\!\!\!\!\!\!\!\!\!\!\!\Gamma_*({\mathcal F}[[\natural]])_{(fg)}$}}%
%
\special{pn 8}%
\special{pa 1000 2190}%
\special{pa 3000 2190}%
\special{fp}%
\special{sh 1}%
\special{pa 3000 2190}%
\special{pa 2934 2170}%
\special{pa 2948 2190}%
\special{pa 2934 2210}%
\special{pa 3000 2190}%
\special{fp}%
\put(31.9000,-22.0000){\makebox(0,0)[lb]{$\Gamma(D_+(fg), {\mathcal F}[[\natural]])$}}%
%
\special{pn 8}%
\special{pa 3390 600}%
\special{pa 3390 2000}%
\special{fp}%
\special{sh 1}%
\special{pa 3390 2000}%
\special{pa 3410 1934}%
\special{pa 3390 1948}%
\special{pa 3370 1934}%
\special{pa 3390 2000}%
\special{fp}%
\end{picture}%

\par\bigskip\noindent
By a similar argment as above, we define an $\proj[[\natural]]$-module homomorphism 
\begin{equation}\label{5.23}
\beta_{\mathcal F}:\tilde{\Gamma_*({\mathcal F}[[\natural]])}\to {\mathcal F}[[\natural]].\end{equation}

\begin{prop}\label{serre}
Assume that a graded ring $S$ is generated by 
$S_1=\{f_1,f_2,\ldots,f_\ell\}~(\exists \ell \in {\mathbf Z}_{\geq 0})$ 
as $S_0$-algebra. Then we see 
\begin{enumerate}
\item[](i) If $S$ is a domain then the map $\alpha$ induced in 
(\ref{5.22}) is injective.
\item[](ii) If $(f_i)_{(i=1,2,\ldots,\ell)}$ are all prime ideals, then 
the map $\alpha$ is an isomorphism. 
\item[](iii) When $S=k[z_0,z_1,z_2,\ldots,z_n]$, then the map 
$\alpha$ is an isomorphism.
\end{enumerate}
\end{prop}

We also have the following: 
\begin{prop}\label{thm5.22}
Assume that ${\mathcal F}[[\natural]]$ 
is a quasi-coherent sheaf $\proj[[\natural]]$-module. 
Then the homomorphism $\beta_{\mathcal F}$ induced 
in (\ref{5.22}) is an isomorphism. 
Furthermore, via (\ref{5.22-1}), we see 
\begin{equation}
\begin{array}{ll}
&H^0({\mathbb P}_{k[z_0,z_1,\ldots,z_n]}^n[[\natural]])\\
:=&H^0({\rm Proj}(k[z_0,z_1,z_2,\ldots,x_n]),
{\mathcal O}_{{\rm Proj}(k[z_0,z_1,z_2,\ldots,x_n])}(m)[[\natural]] ) \\
=&
\left\{
\begin{array}{l}
0\quad\quad\quad\quad\quad\quad\quad (\mbox{if }m<0),  \\
k[z_0,z_1,z_2,\ldots,x_n]|_m \quad (\mbox{o.w.  }m\geq 0).
\end{array}\right.
\end{array}
\end{equation}
\end{prop}

\begin{defn}
Assume that $\Lambda$ is a holomorphic skew-biderivation satisfying Jacobi rule.Then, $\#$ is called {\rm symbol calculus} on ${\rm Proj}(S)$ if 
$({\rm Proj}(S)[[\natural]],\#)$ has an associative algebra sheaf structure 
such that 
$$
f(Z)\# g(Z) 
=f\bullet g+ \left(\frac{\mu}{2}\right)
\Lambda^{\alpha\beta}
\partial_{Z_{\alpha}}f(Z)
\partial_{Z_{\beta}}g(Z)+
\cdots.
$$
\end{defn}

\par\bigskip\noindent
{\bf Proof of Theorem \ref{thm1}}.~
Under the assumption (\ref{lambda-relation}), 
it is easy to check 
\begin{equation}\label{sharp-product-2}
\begin{array}{ll}
&f(Z)\sum_{k=0}^\infty 
\frac{1}{k!}\left(\frac{\mu}{2}\right)^k
\bigr(\stackrel{\leftarrow}
{\partial_{Z_{{\alpha}_1}}}
\Lambda^{{\alpha}_1,{\beta}_1} 
        \stackrel{\rightarrow}{\partial_{Z_{{\beta}_1}}}\bigr)
\cdots
\bigr(\stackrel{\leftarrow}{\partial_{Z_{{\alpha}_k}}}
\Lambda^{{\alpha}_k,{\beta}_k} 
        \stackrel{\rightarrow}{\partial_{Z_{{\beta}_k}}}\bigr)g(Z)\\
=& \sum_{k=0}^\infty \frac{1}{k!}\left(\frac{\mu}{2}\right)^k
\Lambda^{\alpha_1\beta_1}
\Lambda^{\alpha_2\beta_2}
\cdots \Lambda^{\alpha_k\beta_k}
\partial_{Z_{\alpha_1}}\partial_{Z_{\alpha_2}}\cdots \partial_{Z_{\alpha_k}}f(Z)\partial_{Z_{\beta_1}}\partial_{Z_{\beta_2}}\cdots \partial_{Z_{\beta_k}}g(Z). 
\end{array}
\end{equation}
Then the right hand side of (\ref{sharp-product-2}) 
coincides with the asymptotic expansion formula 
for product of the Weyl type pseudo-differential operators. 
Thus, it completes the proof. $\Box$
\par\bigskip\noindent
{\bf Proof of Theorem \ref{thm1-5}}.~
As seen in the previous argument, 
as for sheaf cohomology of projective space, we obtain that   
\begin{equation}\label{sheaf-cohomology}
\sum_{k=0}^\infty H^0(\cpn, {\mathcal O}_\cpn(k))\cong 
\sum_{k=0}^\infty {\mathbb C}[Z]_k, 
\end{equation}
where ${\mathbb C}[Z]_k$ stands for  
the space of homogeneous polynomials of degree $k\in {\mathbf Z}_{\geq 0}$. 
Then a direct computation using (\ref{sharp-product-2}) shows 
that $\mu$ can be 
specialized a scalar. 
$\Box$ 
\par\bigskip\noindent
{\bf Proof of Theorem \ref{thm2}}.~  
We would like to compute exponentials having the following form 
$f(Z)=g(t)e^{\frac{1}{\mu}Q[Z](t)}$ 
with respect to $\#$ for quadratic polynomials 
under a quite general setting.  

Let $Z=[z^1:\ldots:z^{n}]$, $A[Z]:=ZA{}^tZ$, where 
$A \in Sym(n,{\mathbb C})$, i.e. 
$A$ is an $n\times n$-complex symmetric matrix. In order to compute 
the exponential $F(t):=e_{\#}^{t\frac{1}{\mu}A[Z]}$ with respect to the 
Wyel type product formula, 
we treat the following evolution equation: 
\begin{equation}\label{evolution}
\partial_t F=\frac{1}{\mu}A[Z]\# F, 
\end{equation}
with an initial condition 
\begin{equation}\label{initial}
F_0=e^{\frac{1}{\mu}B[Z]},
\end{equation}
where $B\in Sym(n,{\mathbb C})$. 

As seen above, our setting 
is rather different from the situations 
considered in the article \cite{maillard} 
and in \cite{o, ommy1, ommy5}. 
See also \cite{q}\footnote{
Quillen's method employing the Cayley transform is very useful to compute superconnection character forms and 
supertrace of Dirac-Laplacian heat kernels (cf. \cite{gt}).}. 
However, to compute exponentials, 
we can use similar methods employed in the articles above, 
as will be seen below: 

Under the assumption $F(t)=g\cdot e^{\frac{1}{\mu}Q[Z]}$ 
($g=g(t),~Q=Q(t)$), we would like to find a solution of the equations 
(\ref{evolution}) and (\ref{initial}). 

Direct computations give  
{ 
\begin{eqnarray}\label{LHS}
\nonumber
\mbox{L.H.S. of } (\ref{evolution}) &=& 
g' e^{\frac{1}{\mu}Q[Z]} +g {\frac{1}{\mu}Q'[Z]}e^{\frac{1}{\mu}Q[Z]}, \\
\label{RHS}
\nonumber
\mbox{R.H.S. of } (\ref{evolution}) &=& {\frac{1}{\mu}A[Z]}{\#}F  \\
\nonumber
&\stackrel{(\ref{sharp-product})}{=}&
{\frac{1}{\mu}A[Z]}\cdot F+\frac{i\hbar}{2}
\Lambda^{i_1j_1}\partial_{i_1}{\frac{1}{\mu}A[Z]}\cdot \partial_{j_1}F  \\
\nonumber
&&\qquad-\frac{\hbar^2}{2\cdot 4} \Lambda^{i_1j_1}
\Lambda^{i_2j_2} \partial_{i_1i_2}{\frac{1}{\mu}A[Z]}
\partial_{j_1j_2}F\\
\end{eqnarray}}
%
%
where $A=(A_{ij}), \Lambda=(\Lambda^{ij})$ and $Q=(Q_{ij})$.  
Comparing the coefficient of $\mu^{-1}$ gives 
{ 
\begin{equation}
Q'[Z]=A[Z]-2{}^tA\Lambda Q[Z]-Q\Lambda A \Lambda Q[Z] . 
\end{equation}
}
Applying $\Lambda$ by left and setting 
{{} $q:=\Lambda Q$ and $a:=\Lambda A$}, 
we easily obtain 
{{} 
\begin{eqnarray}
\Lambda Q' &=& \Lambda A + \nonumber
\Lambda Q\Lambda A-\Lambda A\Lambda Q 
-\Lambda Q\Lambda A\Lambda Q \\
\nonumber&=&
(1+\Lambda Q)\Lambda A(1-\Lambda Q)\\
&=&(1+q)a(1-q).
\end{eqnarray}
}
As to the coefficient of $\mu^0$, we have 
{{} 
\begin{eqnarray}
\nonumber
g'&=&\frac{1}{2}\Lambda^{i_1j_1}\Lambda^{i_2j_2}A_{i_1i_2}gQ_{j_1j_2}  
\\
&=&-\frac{1}{2}tr(aq) \cdot g ,  
\end{eqnarray}
where $``{\it tr}"$ means the trace.  
Thus we obtain 
\begin{prop}
The equation (\ref{evolution}) is rewritten by 
{{} 
\begin{eqnarray}
\partial_t q &=& (1+q)a(1-q), \label{equation-1}\\ 
\partial_t g &=& -\frac{1}{2}tr(aq)\cdot g. \label{equation-2}
\end{eqnarray}
}
\end{prop}
In order to solve the equations (\ref{equation-1}) and (\ref{equation-2}), 
we  now recall the ``Cayley transform." 
\par\bigskip 
{{} 
\begin{prop}\label{cayley}
Set 
\begin{equation}\label{cayley-transf}
C(X):=\frac{1-X}{1+X}
\end{equation} 
if $\det (1+X)\not=0$ . Then 
\begin{enumerate}
\item $X\in sp_{\Lambda}(n,{\mathbb R})\Longleftrightarrow
\Lambda X\in Sym(n,{\mathbb R})$, \\ and then 
$C(X)\in Sp_{\Lambda}(n,{\mathbb R})$, 
where
\begin{eqnarray}
\nonumber&&Sp_{\Lambda}(n,{\mathbb R})
:=\{g\in GL(n,{\mathbb R})|{}^tg\Lambda g=\Lambda\},\\
\nonumber&&sp_{\Lambda}(n,{\mathbb R})
:=Lie(Sp_{\Lambda}(n,{\mathbb R})).
\end{eqnarray}
\item $C^{-1}(g)=\frac{1-g}{1+g}$, (the ``{\rm inverse Cayley transform}").
\item $e^{2\sqrt{-1}a}=c(-\sqrt{-1}\tan(a))$. 
\item $\log a = 2\sqrt{-1}\arctan (\sqrt{-1}C^{-1}(g))$. 
\item $\partial_t q =(1+q)a(1-q)\Longleftrightarrow 
\partial_t C(q)=-2aC(q).$ 
\end{enumerate}
\end{prop}
}

Solving the above equation 5 in Proposition \ref{cayley}, 
we have 
{{} $$C(q)=e^{-2at}C(b),$$} 
\noindent
where $b=\Lambda B$ and then 
$$
q=C^{-1}\bigr(e^{-2at}\cdot C(b)\bigr)
=C^{-1}\bigr(C(-\sqrt{-1}\tan(\sqrt{-1}at))\cdot C(b )\bigr).
$$
Hence, according to the inverse Cayley transform, 
we can get $Q$ in the following way. 
\begin{prop}
\begin{equation}\label{Q}
Q=-\Lambda\cdot C^{-1}\Bigr(C(-\sqrt{-1}\tan (\sqrt{-1} \Lambda A t)) 
\cdot C(\Lambda B)\Bigr). 
\end{equation} 
\end{prop}
Next we compute the amplitude coefficient part $g$. 
Solving 
\begin{equation}
g'=-\frac{1}{2}tr(aq)\cdot g
\end{equation}
gives 
\begin{prop}
\begin{equation}
g=\det{}^{-\frac{1}{2}}
\Bigr(\frac{e^{at}(1+b)+e^{-at}(1-b)}{2}\Bigr). 
\end{equation} 
\end{prop}
Setting $t=1$, $a=\Lambda A$ and $b=0$, 
we get 
\begin{thm}\label{star-exponent}
{ 
\begin{eqnarray}
\label{star_exponential}
e_{\#}^{\frac{1}{\mu}A[Z]}
&=&\det{}^{-\frac{1}{2}}\Bigr(\frac{e^{\Lambda A}+e^{-\Lambda A}}{2}\Bigr)
\cdot e^{\frac{1}{\mu}
(\frac{\Lambda^{-1}}{\sqrt{-1}} \tan (\sqrt{-1}\Lambda A )) [Z]  } .
\end{eqnarray}
}
\end{thm}
As usual, using the $\check{\rm C}$ech resolution, 
we can compute the sheaf-cohomology 
$\sum_{k=0}^\infty H^0(\cpn, {\mathcal O}_\cpn[\mu,\mu^{-1}]])$. 
Combining it with Theorems 
\ref{star-exponent}, 
we see 
$$e_{\#}^{\frac{1}{\mu}A[Z]} \in 
\sum_{k=0}^\infty H^0(\cpn, {\mathcal O}_\cpn[\mu,\mu^{-1}]]). $$   
This completes the proof of Theorem~\ref{thm2}. $\Box$
\par\bigskip\noindent

\section{Remarks} 
{
Consider the following diagram: 
\par\bigskip
\qquad\qquad\quad
\unitlength 0.1in
\begin{picture}( 30.0000, 23.1500)( 16.0000,-26.3000)
%
\special{pn 8}%
\special{pa 3200 800}%
\special{pa 1800 2400}%
\special{fp}%
\special{sh 1}%
\special{pa 1800 2400}%
\special{pa 1860 2364}%
\special{pa 1836 2360}%
\special{pa 1830 2338}%
\special{pa 1800 2400}%
\special{fp}%
%
\special{pn 8}%
\special{pa 3200 800}%
\special{pa 4600 2400}%
\special{fp}%
\special{sh 1}%
\special{pa 4600 2400}%
\special{pa 4572 2338}%
\special{pa 4566 2360}%
\special{pa 4542 2364}%
\special{pa 4600 2400}%
\special{fp}%
\put(32.0000,-4.0000)
{\makebox(0,0){$ \bigr((x^{\alpha,\dot{\alpha}}),[\pi_1:\pi_2]\bigr)
\in M:={\mathbb C}^{4}\times{\mathbb C}{\mathbb P}^1$}}%
\put(16.0000,-28.0000){\makebox(0,0)[lb]{$([z_1:\ldots:z_4])\in\sts$}}%
\put(44.0000,-28.0000){\makebox(0,0)[lb]{$(x^{\alpha,\dot{\alpha}})
\in{\mathbb C}^{4}$}}%
\put(44.0000,-14.0000){\makebox(0,0)[lb]{$\Pi_2$}}%
\put(16.0000,-14.0000){\makebox(0,0)[lb]{$\Pi_1$}}%
\end{picture}%

\par\vspace{1cm}
\noindent
where $x^{\alpha,\dot{\alpha}}$ are even variables, 
we set 
\begin{eqnarray}
&&\nonumber (x^{\alpha,\dot{\alpha}})
:=(x^{1,\dot{1}},x^{1,\dot{2}},x^{2,\dot{1}},x^{2,\dot{2}}),\\
&&\nonumber ([z_1:\ldots:z_4]) 
:=([x^{\alpha,\dot{1}}\pi_\alpha:x^{\alpha,\dot{2}}\pi_\alpha:\pi_1:\pi_2]). 
\end{eqnarray}
Here we use Einstein's convention
(we will often omit $\sum$ unless there is a danger of confusion). 
We call 
$([z_1:\ldots:z_4])$ 
the {\rm homogeneous coordinate system} of $\sts$. 
\begin{enumerate}
\item 
The relations\footnote{Here $[~,~]$ denotes the commutator bracket.} 
($\dot{\alpha}, \dot{\beta} = \dot{1}, \dot{2}$)  
\begin{equation}
[z^{\dot{\alpha}},z^{\dot{\beta}}]=
\hbar D^{\alpha\dot{\alpha},\beta\dot{\beta}}\pi_\alpha\pi_\beta, 
\end{equation}
where $z^{\dot{1}}:=z_1,~z^{\dot{2}}:=z_2$, 
give a globally defined non-commutative associative product
$\#$ on ${\mathbb C}{\mathbb P}^3$, 
where $\bigr(D^{\alpha\dot{\alpha},\beta\dot{\beta}}\bigr)$ 
is a skew symmetric matrix. 
\item 
Let $A[Z]$ be a homogeneous polynomial
of $z^{\dot{1}}=z_1=x^{\alpha,\dot{1}}\pi_{\alpha},~
z^{\dot{2}}=z_2=x^{\alpha,\dot{2}}\pi_{\alpha}$ 
with degree $2$. 
Then a star exponential function $e_{\#}^{\frac{1}{\mu}A[Z]}$ 
gives a ``function"
on ${\mathbb C}{\mathbb P}^{3}$. 
\end{enumerate} 
}
More precisely, 
\begin{thm}
Assume that $\Lambda:=\hat{\Lambda}$ 
and $A[Z]$ a homogeneous polynomial of 
$z^{\dot{1}}=x^{\alpha,\dot{1}}\pi_{\alpha},~
z^{\dot{2}}=x^{\alpha,\dot{2}}\pi_{\alpha}$ with degree $2$. 
Then a star exponential function $e_{\#}^{\frac{1}{\mu}A[Z]}$ 
gives a cohomology class of 
${\mathbb C}{\mathbb P}^3$ with coefficients in the sheaf 
$\sum_{k=0}^{\infty}{\cal O}_{{\mathbb C}{\mathbb P}^3}(k)$. 
\end{thm}

\par\bigskip


\unitlength 0.1in
\begin{picture}( 36.8500, 22.5000)(  6.1500,-31.6500)
\put(30.0000,-10.0000){\makebox(0,0){$M:={\mathbb C}^{4|4N}\times{\mathbb C}{\mathbb P}^1$}}%
%
\special{pn 8}%
\special{pa 3000 1200}%
\special{pa 3000 1800}%
\special{fp}%
\special{sh 1}%
\special{pa 3000 1800}%
\special{pa 3020 1734}%
\special{pa 3000 1748}%
\special{pa 2980 1734}%
\special{pa 3000 1800}%
\special{fp}%
\put(30.0000,-21.0000){\makebox(0,0){${\mathbb C}^{4|2N}\times{\mathbb C}{\mathbb P}^1$}}%
%
\special{pn 8}%
\special{pa 2800 2400}%
\special{pa 1800 3000}%
\special{fp}%
\special{sh 1}%
\special{pa 1800 3000}%
\special{pa 1868 2984}%
\special{pa 1846 2974}%
\special{pa 1848 2950}%
\special{pa 1800 3000}%
\special{fp}%
%
\special{pn 8}%
\special{pa 3300 2400}%
\special{pa 4300 3000}%
\special{fp}%
\special{sh 1}%
\special{pa 4300 3000}%
\special{pa 4254 2950}%
\special{pa 4254 2974}%
\special{pa 4234 2984}%
\special{pa 4300 3000}%
\special{fp}%
\put(18.0000,-32.5000){\makebox(0,0){${\mathbb C}{\mathbb P}^{3|N}$}}%
\put(43.0000,-32.5000){\makebox(0,0){${\mathbb C}^{4|2N}$}}%
\put(33.5000,-15.0000){\makebox(0,0){$\Pi$}}%
\put(16.0000,-27.0000){\makebox(0,0){$\Pi_1$}}%
\put(43.0000,-27.0000){\makebox(0,0){$\Pi_2$}}%
\end{picture}%

\par\bigskip
\noindent
where $\Pi$ denotes the chiral projection,  
we can consider non-anti-commutative deformation of super twistor space. 

In order to give a brief explanation, 
we recall the definition of super 
twistor manifold (\cite{lebr, taniguchi-miyazaki, ward-wells}). 
\begin{defn}
$(3\vert N)$-dimensional complex super manifold $Z$ is said to be 
a super twistor space if the following conditions $(1)-(3)$ are satisfied.  
\begin{enumerate}
\item[]$(1)$ $p: Z\longrightarrow {\mathbb C}{\mathbb P}^1$ 
is a holomorphic fiber bundle.
\item[]$(2)$ $Z$ has a family of holomorphic section of $p$ 
whose normal bundle is isomorphic to 
${\mathcal O}_{ {\mathbb C}{\mathbb P}^1}(1)
\oplus {\mathcal O}_{{\mathbb C}{\mathbb P}^1}(1)\oplus C^N \otimes 
\Pi {\mathcal O}_{{\mathbb C}{\mathbb P}^1}(1)$. 
\item[]$(3)$ $Z$ has an anti-holomorphic involution $\sigma$ being 
compatible with $(1),~(2)$ and $\sigma$ has no fixed point.  
\end{enumerate}
\end{defn}


We define 
${\mathbb C}{\mathbb P}^{3|N}_{*_{\alpha'}}=({\mathbb CP}^{3|N}, 
{\mathcal O}_{{\mathbb C}{\mathbb P}^{3|N}, *_\alpha'})$. 

Let $f(z\vert \xi ; \alpha^{\prime})$ be 
a local section defined in the following manner: 
\begin{equation}\label{local-represent-2}
 f(z\vert \xi ; \alpha^{\prime})=\sum_{k=0}^{N} \sum_{1\leq i_1 \leq i_2 \leq 
 \ldots \leq 
i_k \leq N} f_{i_1 i_2 \ldots i_k} (z) 
\xi^{i_1}  \xi^{i_2}  \cdots  \xi^{i_k}  
\end{equation}
where $f_{i_1\ldots i_k}(z)$ is 
a homogeneous element of $z=[z_1:z_2:z_3:z_4]$ with homogeneous degree 
$(-k)$ on ${\mathbb C}{\mathbb P}^3$. 
Then we can introduce a structure 
sheaf ${\mathcal O}_{{\mathbb C}{\mathbb P}^{3\vert N}, *\alpha^{\prime}}$
whose local section is given by $f(z\vert \xi ; \alpha^{\prime})$.  

Under these notations, we can introduce a ringed space denoted by 
${\mathbb C}{\mathbb P}^{3\vert N}_{*\alpha^{\prime}} = 
({\mathbb C}{\mathbb P}^3, 
{\mathcal O}_{{\mathbb C}{\mathbb P}^{3\vert N}, *_{\alpha'}})$. 
We shall call it non-anti-commutative complex projective super space. 

As for the non-anti-commutative deformed product $*$ 
aassociated with the non-anti-commutative complex projective super space,  
we have commutation relations of local coordinate functions: 
\begin{enumerate}
\item[]$(1)$ Let $ (z_1, z_2, \pi_1, \pi_2 \vert \xi^1,\ldots,\xi^N )$ 
be a local coordinate system of ${\mathcal P}^{3\vert N}$, 
where ${\mathcal P}^{3\vert N}$ denotes the non-anti-commutative 
open super twistor space. 
Then 
\[ \{ \xi^i ,\xi^j \}_* 
=\alpha^{\prime} C^{i\alpha, j\beta}\pi_{\alpha}\pi_{\beta}, \qquad 
(0~ \rm{o.w.})\]
\item[]$(2)$ A local coordinate system $(z_1, z_2, \pi_1, \pi_2 
\vert \xi^1,\ldots, \xi^N )$ 
of ${\mathbb C}{\mathbb P}^{3\vert N}$ satisfies 
\[ \{ \xi^i ,\xi^j \}_* 
=\alpha^{\prime} C^{i\alpha, j\beta}\pi_{\alpha}\pi_{\beta}, \qquad 
(0~\rm{o.w.}) \]
\end{enumerate}

Here we do not explain more the notion and notations which appeared above 
and do not give the proof of them. 
For details, see \cite{taniguchi-miyazaki}.


\begin{thebibliography}{99}

\bibitem{bffls} F. Bayen, M. Flato, C. Fronsdal, A. Lichnerowicz, 
and D. Sternheimer, {\it Deformation theory and quantization I,}
Ann. of Phys. 111 (1978), 61-110.



\bibitem{fedosov}B. V. Fedosov, 
{\it A simple geometrical construction of deformation quantization, }
Jour. Diff. Geom. 40 (1994), 213-238.

\bibitem{gt} K. Gomi and Y. Terashima, 
{\em Chern-Weil construction for twisted K-theory}, 
preprint.

\bibitem{homma} Y. Homma and T. Tate, private communications. 

\bibitem{hormander-fio} L. H\"{o}rmander, 
{\em Fourier integral operators. I,}
Acta Mathematica 127, (1971), 79-183.

\bibitem{hormander} L. H\"{o}rmander, 
{\em The Weyl culculus of general pseudodifferential operators,}
Commun. Pure. Appl. Math. 32 (1979), 359-443. 

\bibitem{hormander-I} L. H\"{o}rmander, 
{\em The Analysis of Linear Partial Differential Operators I,}
Springer-Verlag, 1989 

\bibitem{hormander-III} L. H\"{o}rmander, 
{\em The Analysis of Linear Partial Differential Operators III,}
Springer-Verlag, 1984 




\bibitem{k}M. Kontsevich, 
{\it Deformation quantization of Poisson manifolds,} 
Lett. Math. Phys. 66, 
(2003), 157-216. 

\bibitem{lebr}C. LeBrun, Y. S. Poon and R. O. Wells, Jr.,
{\em Projective embedding of complex supermanifolds,}
Commun. Math. Phys. 126 (1990) 433-452.


\bibitem{ommy0} Y. Maeda, N. Miyazaki, H. Omori and A. Yoshioka,
   {\it Star exponential functions as two-valued elements}, 
   {Progr. Math. 232} (2005), 483-492, Birkh\"{a}user. 

\bibitem{maillard} J. M. Maillard, 
{\em Star exponential for any ordering of the elements of 
the inhomogeneous symplectic Lie algebra}, 
Jour. Math. Phys. vol.45, no. 2, (2004), 785-794. 





\bibitem{miyazaki0}N. Miyazaki, 
{\em On a certain class of oscillatory 
integral transformations which determine canonical graphs}, 
Japanese Jour. Math. 24, no.1 (1998) 61-81. 

\bibitem{miyazaki00}N. Miyazaki, 
{\em Formal deformation quantization and the Index theorem}, 
unpublished, (2003), 
in Japanese. 


\bibitem{miyazaki01}N. Miyazaki, 
{\em Lifts of symplectic diffeomorphisms 
of a Weyl algebra bundle with Fedosov connection}, 
International Journal of Geometric Method in Modern Physics,  
Vol. 4, No. 4 (2007) 533-546.


\bibitem{miy}N. Miyazaki, 
{\it A Lie Group Structure for Automorphisms of 
a Contact Weyl Manifold}, Progr. Math. 252 (2007), 
25-44, Birkh\"{a}user. 


\bibitem{miyazaki-rims-2009} N. Miyazaki, {\em 
Remarks on deformation quantization}, Kyoto University RIMS Kokyuroku 1692, 
Geometric Mechanics, (2010), 1-16. 





\bibitem{miyazaki-preparation}N. Miyazaki, 
{\em Deformation of structure sheaf on a projective shceme,} 
in preparation. 


\bibitem{moriyoshi} H. Moriyoshi, private communications. 

\bibitem{murphy} G. J. Murphy, {\em $C^*$-algebras and operator theory}, 
Academic Press, (1990)

\bibitem{o} H. Omori, {\em Physics in Mathematics}, The Univ. Tokyo Press, 
(2004), in Japanese. 



\bibitem{ommy1} H. Omori, Y. Maeda, N. Miyazaki and A. Yoshioka, 
   {\it Star exponential functions for quadratic forms and polar elements}, 
   {Contemporary Mathematics. 315 (Amer. Math. Soc.)}, (2002), 25-38.   








\bibitem{ommy4} H. Omori, Y. Maeda, N. Miyazaki and A. Yoshioka, 
{\em Convergent star product 
on Fr{\'e}chet-Poisson algebras of Heisenberg type,}  
Contemporary Mathematics 434 (Amer. Math. Soc.), (2007), 99-123.   


   
\bibitem{ommy5} H. Omori, Y. Maeda, N. Miyazaki and A. Yoshioka, 
{\it A new nonformal noncommutative calculus:
Associativity and finite part regularization}, 
{Ast\'{e}risque} 321 (2008), 267-297. 



\bibitem{q} D. Quillen, 
{\em Superconnection character forms and the Cayley transform}, 
Topology, 27 (2), (1988), 211-238. 

\bibitem{sato} H. Sato, private communications. 


\bibitem{shubin} M. A. Shubin, {\em Pseudodifferential 
Operators and Spectral Theory,}Springer-Verlag. 

\bibitem{s} D. Sternheimer,  
{\it Deformation quantization twenty years after,}
AIP Conf. Proc. 453 (1998), 107-145.


\bibitem{suzuki} T. Suzuki, private communications. 

\bibitem{taniguchi-miyazaki} T. Taniguchi and N. Miyazaki, 
{\em On non(anti)commutative super twistor spaces}, 
International Journal of Geometric Method in Modern Physics 
Vol. 7, No. 4 (2010) 655-668. 

\bibitem{ward-wells}R. S. Ward and R. O. Wells, Jr.,
{\em Twistor geometry and field theory,}
Cambridge Monographs on Math. Phys. 1990.

\bibitem{w} N. Woodhouse,  {\it Geometric quantization,} Clarendon Press,  
(1980), Oxford.



\bibitem{y} A. Yoshioka,  
   {\it Contact Weyl manifold over  a symplectic manifold,} 
    in ``Lie groups, Geometric structures and 
    Differential equations", Adv. Stud. Pure Math. 37(2002), 459-493. 


\end{thebibliography}
\end{document}